\documentclass[11pt,a4,twoside]{article}
\usepackage{amsmath}
\usepackage{amsfonts}
\usepackage{amssymb}
\usepackage[top=1in, bottom=1in, left=1in, right=1in]{geometry}
\usepackage{graphicx}
\usepackage[bottom]{footmisc}
\usepackage{fancyhdr}
\begin{document}
	\begin{titlepage}
		\begin{center}
		\huge{\bfseries The asymptotic properties of $\phi(n)$ and a problem related to visibility of Lattice points}\\
		[1.5 cm]
		\Large{Debmalya Basak\footnote{Email id : db15ms149@iiserkol.ac.in}}\\
		[4mm]
		\Large{Indian Institute of Science Education and Research,Kolkata}\\
		[0.75cm]
		\vspace{5cm}
        \begin{abstract}
        	We look at the average sum of the Euler's phi function $\phi{(n)}$  and it's relation with the visibility of a point from the origin.We show that $\forall{\hspace{0.05in}{k} \ge{1}} ,k\in\mathbb{N},\exists$ a \textit{k$\times$k} grid in the 2D space such that no point inside it is visible from the origin.We define visibility of a lattice point from a set and try to find a bound for the cardinality of the smallest set S such that for a given \textit{n} $\in\mathbb{N}$,all points from the \textit{n$\times$n} grid are visible from S.  
        \end{abstract}
\pagenumbering{roman}
		\end{center}
	\end{titlepage}
\newpage
\tableofcontents
\newpage
\pagenumbering{arabic}
\section{On the Visibility of Lattice Points in 2-D Space}
\lhead{\textit{VISIBILITY OF LATTICE POINTS IN 2 D SPACE}}
\rhead{\textit{Debmalya Basak}}
\hspace{0.6cm}First let's define the Euler totient function,something which will be very useful in this section.If $n \geq 1$,we define $\phi {(n)}$ as the number of positive integers less than $n$ and coprime to $n$.Now,we will introduce the concept of visibility of lattice points.We say that 2 integer lattice points ($a,b$) and ($c,d$) are visible if the line joining those 2 points doesn't contain any other lattice point in between.Now,we prove a very important result here.
\vspace{3mm}

\hspace{1.2cm}\textbf{Theorem 1 :} 2 integer lattice points ($a,b$) and ($c,d$) are visible if and only if ($a-c$) 

\hspace{3.6cm} and ($b-d$) are coprime to each other.

\vspace{3mm}\hspace{2.35cm}\textit{Proof }: The proof of this theorem is very easy.To make the problem simplier,we take ($c,d$) to be the origin as we can translate our 2 D space to our wish.We first prove the $\implies$ side i.e if ($a,b$) is visible from the origin,then ($a,b$)$=1$.Suppose $a$ and $b$ are not coprime to each other.Then\hspace{1mm} $\exists$\hspace{1mm} $k>1$,$k$ $\in\mathbb{N}$ such that $k|a$ and $k|b$.Let $a=ka_1$ and $b=kb_1$.Then the lattice point ($a_1,b_1$) satisfies the equation of the line joining ($a.b$) and the origin and lies in between the 2 endpoints.So,we arrive at a contradiction.

Now,the $\impliedby$ side i.e. if ($a,b$)$=1$,then ($a,b$) is visible from the origin.It goes the same way by assuming that ($a,b$) is not visible and arriving at a trivial contradiction.We skip this part and leave it to the reader.

Now that we get a feeling that the visibility of lattice points is highly involved with the Euler totient function,let's study $\phi{(n)}$ a little more in details.
\subsection{Average order of the Euler Totient function}
\hspace{0.6cm}As in most cases,arithmetic functions are very irregular of the domain of natural numbers.So,it is always better to study the average of such funcions over a certain range which gives us a nice picture of certain properties of the function.Same is the case with $\phi{(n)}$.In this part,we will try to estimate $\sum \limits_{\substack{ n\leq x}} \phi{(n)}$ for some $x$ $\in\mathbb{R}$.We have some trivial properties of $\phi{(n)}$ which are very important for mathematical purpose.We will not state the proofs here.

\vspace{3mm}\hspace{0.5cm}\textbf{1.} $\sum \limits_{\substack{ d|n}} \phi{(d)}=n$ \hspace{0.5cm}\textbf{2.} $\phi{(n)}=\sum \limits_{\substack{ d|n}} \mu{(d)}$($\dfrac{n}{d}$)\footnote{$\mu{(d)}$ is known as the M\H{o}bius function and is defined as follows.For any natural number $n$,$\mu{(n)}=1$,if $n$ is squarefree and has an even number of prime factors,$\mu{(n)}=-1$ if $n$ is squarefree and has an odd number of prime factors and finally,$\mu{(n)}=0$ if $n$ is not squarefree.}  \hspace{0.5cm} \textbf{3.} $\phi{(n)}=\phi{(a)}\phi{(b)}$ if $n=ab$ and ($a,b$)$=1$
\vspace{3mm}

And,finally,we give the general formula of $\phi{(n)}$.

\vspace{3mm}
\hspace{1cm}$\phi{(n)}=n\prod_{p|n}$(1$-\dfrac{1}{p}$),where $p$ is a prime,$n$ $\in\mathbb{N}$.
\vspace{3mm}

Now,we can easily write that  $|\mu{(n)}| \leq 1$ \hspace{1mm}$\forall \hspace{1mm} n \in\mathbb{N}$.So,we can say that the $\sum_{n=1}^{\infty} |\dfrac{\mu{(n)}}{n^2}|$ is absolutely convergent\footnote{We assume that we know this.For more,look at References.} as $\sum_{n=1}^{\infty} \dfrac{1}{n^2}=\dfrac{\pi^2}{6}=\zeta{(2)}$\footnote{Here,$\zeta{(n)}$ is the Reimann Zeta function.}.So,the sum  $\sum_{n=1}^{\infty} \dfrac{\mu{(n)}}{n^2}$ is also convergent.We use the value of this summation  $\sum_{n=1}^{\infty} \dfrac{\mu{(n)}}{n^2}=\dfrac{6}{\pi^2}$ for the estimation of the average order of $\phi{(n)}$,the proof of which we won't state here.

We have\footnote{We assume here,that $O(\sum \limits_{\substack{ n>x}}\dfrac{1}{n^2})=O(\dfrac{1}{x})$},

\vspace{3mm}
\hspace{1cm}$\sum \limits_{\substack{ n \leq x}}\dfrac{\mu{(n)}}{n^2} = \sum_{n=1}^{\infty} \dfrac{\mu{(n)}}{n^2}-\sum \limits_{\substack{ n>x}}\dfrac{\mu{(n)}}{n^2}= \dfrac{6}{\pi^2}+O(\sum \limits_{\substack{ n>x}}\dfrac{1}{n^2})=\dfrac{6}{\pi^2}+O(\dfrac{1}{x})$

\vspace{3mm}

So,now,we get\footnote{By Euler Summation method,we can prove that $\sum \limits_{\substack{ d \leq x}}\dfrac{1}{d}=\log x +C+O$($\dfrac{1}{x}$),where $C>0$ is a constant.\\

	\hspace{0.35cm}Also,$\sum \limits_{\substack{ q \leq x}}q= \dfrac{x^2}{2}+O$($x$).} ,

\vspace{3mm}\hspace{1cm}$\sum \limits_{\substack{ n\leq x}} \phi{(n)}=\sum\limits_{\substack{n< x}}$$\sum\limits_{\substack{d|n}}\mu{(d)}$ ($\dfrac{n}{d})=\sum \limits_{\substack{ q,d\\ qd \leq x}} \mu{(d)}q=\sum \limits_{\substack{ d\leq x}} \mu{(d)}\sum \limits_{\substack{ q\leq x/d}} q=\sum \limits_{\substack{ d\leq x}} \mu{(d)}$\{$\dfrac{1}{2}$($\dfrac{x}{d}$)$^2$$+O$($\dfrac{x}{d}$)\}

\vspace{3mm}\hspace{2.5cm}$=\dfrac{x^2}{2}\sum \limits_{\substack{ d \leq x}}\dfrac{\mu{(d)}}{n^2}+O$($x\sum \limits_{\substack{ d \leq x}}\dfrac{1}{d}$)$=\dfrac{3x^2}{\pi^2}+O$($x\log x$)

\vspace{3mm}
Hence,we can say that the average order of $\phi{(n)}$ is $\dfrac{3n}{\pi^2}$.

\subsection{Density of Lattice points visible from the origin}
\rhead{\textit{Debmalya Basak}}
\lhead{\textit{DENSITY OF LATTICE POINTS VISIBLE FROM THE ORIGIN}}
\hspace{0.6cm}Now that we are a bit familiar with the average order of the Euler Totient function,we look at one of it's most important applications.Suppose,we look at an $n\times n$ grid.The total number of points in this grid is $n^2$.Let us define a set $A_n$ which contains all the points lying in this $n\times n$ grid which are visible from the origin.We have a very interesting theorem regarding this set $A_n$.

\vspace{3mm}
\hspace{1.2cm}\textbf{Theorem 2 :} The set of all integer lattice points visible from the origin has asymptotic 

\hspace{3.75cm}density $\dfrac{6}{\pi^2}$ i.e. $\lim_{n \to \infty} \dfrac{|A_n|}{n^2}=\dfrac{6}{\pi^2}$.

\vspace{3mm}\hspace{2.35cm}\textit{Proof }: We make a few observations at the beginning.All the 8 lattice points adjacent to the origin are visible from the origin.Now,we divide the 2 D space into 8 subdivisions along the straight lines $y=x$,$y=-x$,$y=0$ and $x=0$.By symmetry,we can look at the number of lattice points visible from the origin in any one of these 8 subdivisionsand multiply the count by 8.Suppose we look at the subdivision lying in the region marked by the 2 straight lines $y=0$ and $y=x$ except the 2 points adjacent to the origin i.e. ($1,1$) and ($1,0$).Let us now fix $k>0,k \in\mathbb{N}$. So,we are now looking for integer lattice points ($a,b$) such that $2 \leq a \leq k,1 \leq b<a$ and ($a,b$)$=1$.This is because any point on the 2 boundary lines other than ($1,1$) and ($1,0$) is invisible from the origin.This total count,including the adjacent point to the origin is nothing but $\sum \limits_{\substack{ 1 \leq n\leq k}} \phi{(n)}$ for each subdivision.So,in total,the number of lattice points in this $2k\times2k$ grid is $8$ $\sum \limits_{\substack{ 1 \leq n\leq k}} \phi{(n)}$ while the total number of points is $4k^2$.Hence the asymptotic density of the set$=\lim_{k \to \infty} \dfrac{8\sum \limits_{\substack{ 1 \leq n\leq k}} \phi{(n)}}{4k^2}=\lim_{k \to \infty}\dfrac{\dfrac{24k^2}{\pi^2}+O(k\log k)}{4k^2}=\dfrac{6}{\pi^2}+\lim_{k \to \infty}\dfrac{1}{4}O(\dfrac{\log k}{k})=\dfrac{6}{\pi^2}$\hspace{6.2cm}$\Box$
\vspace{3mm}

So,the probability that a random integer lattice point chosen from the 2 D space is visible from the origin is $\dfrac{6}{\pi^2}$.We can also say that if a random lattice point ($a,b$) is chosen from the 2 D space,the probability that ($|a|,|b|$)$=1$ is $\dfrac{6}{\pi^2}=\dfrac{1}{\zeta{(2)}}$.We can extend this concept from 2 D space to $n$ dimensional space.Interestingly,in $n$ dimensional space,the probability of finding a lattice point visible from the origin is $\dfrac{1}{\zeta{(n)}}$.We will look into a more interesting problem regarding the visibility of lattice points in the next portion.
\section{Hidden Trees in the Forest}
\rhead{\textit{Debmalya Basak}}
\lhead{\textit{HIDDEN TREES IN THE FOREST}}
\hspace{0.6cm}This is a very interesting problem with a lot of research still going on regarding a few aspects of this problem.We also need the concept of Chinese Remainder Theorem for this problem,which we will only state here.
\vspace{3mm}

\textbf{Chinese Remainder Theorem :} Let $m_1,m_2,m_3...m_k$ are positive integers relatively prime to each other i.e. ($m_i,m_j$)$=1$ if $i\neq j$.Suppose $M=m_1m_2m_3...m_k$.Let $a_1,a_2,a_3..a_k$ be a set of arbitrary integers such that $x\equiv a_1\mod m_1$,$x\equiv a_2\mod m_2$,$x\equiv a_3\mod m_3$...$x\equiv a_k\mod m_k$ form a system of linear congruences.Then,the system has an unique solution modulo $M$.

\vspace{3mm}Now,we state the theorem regarding the visibility of lattice points.

\vspace{3mm}
\hspace{1.2cm}\textbf{Theorem 3 :} Given any $k>0,k \in\mathbb{N}$,\hspace{1mm}$\exists$\hspace{1mm} a $k\times k$ grid in the 2 D space such that no 

\hspace{3.75cm}lattice point from that grid is visible from the origin.

\vspace{3mm}
The theorem states that there are large regions in the 2 D space which are completely invisible from the origin.We call these spaces as dark spaces or hidden trees.We now look at a very elegant proof of this,using Chinese Remainder Theorem.

\vspace{3mm}\hspace{2.35cm}\textit{Proof }: It is sufficient to show that for any $k>0,k \in\mathbb{N}$,\hspace{1mm}$\exists$\hspace{1mm} a lattice point ($a,b$) such that none of the lattice points ($a+r,b+s$),$0<r\leq k,0<s\leq k$ is visible from the origin.We construct a $k\times k$ matrix whose first row contains the first $k$ primes,the 2nd row contains the next $k$ primes and so on.Let $d_i$ be the product of all the elements of the $i$th row and $D_i$ be the product of all the elements of the $i$th column.Clearly,($d_i,d_j$)$=1$ and ($D_I,D_J$)$=1$ when $i\neq j$.

Now,we look at the set of congruences $x\equiv -1\mod d_1,x\equiv -2\mod d_2,x\equiv -3\mod d_3....x\equiv -k\mod d_k$.Since $d_i$'s are relatively prime to each other,there exists an unique solution $a$ modulo $d_1d_2..d_k$.We do the same thing for the $D_i$'s.So,we look at the set of congurences $y\equiv -1\mod D_1,y\equiv -2\mod D_2,y\equiv -3\mod D_3....y\equiv -k\mod D_k$ and we have an unique solution $b$ modulo $D_1D_2D_3...D_k$ to the system.Also,we observe that $d_1d_2..d_k=D_1D_2D_3...D_k$.Now,we show that none of the lattice points ($a+r,b+s$),$0<r\leq k,0<s\leq k$ is visible from the origin.Suppose,\hspace{1mm}$\exists$\hspace{1mm} a lattice point ($a+p,b+q$),$0<p\leq k,0<q\leq k$ such that it is visible from the origin.Then,$a\equiv -p\mod d_i$ and $b\equiv -q\mod D_j$ for some $0<i,j \leq k$.But $d_i$ and $D_j$ have a common prime factor $p_{ij}$ because each row and each column has at least one common entry.Hence,($a+p,b+q$)$\geq p_{ij}$.So,we arrive at a contradiction.So,there are arbitrarily large spaces which are completely invisible from the origin.\hspace{15.5cm}$\Box$

We would now move on to more interesting topics regarding this concept.
\section{More interesting problems about the visibility of lattice points}
\hspace{0.6cm}One question asked regarding this concept is about finding a set of lattice points from where the whole $n\times n$ grid is visible.Let's frame the problem mathematically.We say that a set $P$ is completely visible from a set $Q$ if for any lattice point ($a,b$) in $P$,\hspace{1mm}$\exists$\hspace{1mm} a lattice point ($c,d$) in $Q$ such that ($a,b$) is visible from ($c,d$).Suppose we define a set $A_n$ as the set containing all integer lattice points within the $n\times n$ grid.Let $B_n$ be the subset of $A_n$ with least cardinality such that the set $A_n$ is completely visible from $B_n$.We try to study the cardinality of $B_n$ and how it varies with $n$ and we get some interesting results.
\subsection{Abbott's Theorem}
\rhead{\textit{Debmalya Basak}}
\lhead{\textit{ABBOTT'S THEOREM}}

\hspace{0.6cm}H.L. Abbott gave a very strong result regarding the cardinality of the set $B_n$.

\vspace{3mm}
\hspace{1.2cm}\textbf{Theorem 4 :} Let the cardinality of the set $B_n$ be $f$($n$).Then

\vspace{3mm}
\hspace{4cm}
$\dfrac{\log n}{2\log \log n}<f$($n$)$<4\log n$,for sufficiently large $n$.

\vspace{3mm}
We will only provide the proof of the lower bound here and give a sketch of the proof of the upper bound.It follows pretty much the same approach as the proof of the hidden trees.

\vspace{3mm}\hspace{2.35cm}\textit{Proof }: Let $r= \lfloor\dfrac{\log n}{2\log \log n}\rfloor$$+1$,where $\lfloor x \rfloor$ denote the floor function.\\Let ($a_1,b_1$),($a_2,b_2$),($a_3,b_3$)...($a_r,b_r$) be $r$ distinct points from the set $A_n$.We will show that there exists a lattice point in the set $A_n$ which is not visible from any of these $r$ distinct points.Let $p_1,p_2...p_r$ be the first $r$ primes.We look at the following 2 systems of linear congruences,$x\equiv a_i\mod p_i$ and $y\equiv b_i\mod p_i$ \hspace{1mm}$\forall$\hspace{1mm}$1 \leq i \leq r$.By Chinese Remainder Theorem,the 2 systems give us 2 unique solutions $x_0$ and $y_0$ such that ($x_0,y_0$) is different from each of the $r$ points taken initially.We can ensure this disticntion only if $0\leq x_0,y_0 \leq (r+1)p_1p_2..p_r$.Becaause then,we have the oppurtunity to select at least $r+1$ points,one from each interval [$kp_1p_2..p_r,(k+1)p_1p_2..p_r$],$0 \leq k \leq r,k \in\mathbb{N}$ which solves both the systems of linear congruences.But we have taken at most $r$ distinct points initially.So,by Pigeon Hole Principle,at least one of them will be different form all of them.Also,just like Theorem 5,we can show that this new point ($x_0,y_0$) is not visible from any of the $r$ distinct points taken at the beginning.So,we are left to show that ($x_0,y_0$) belongs to the set $A_n$.This is trivial.We use the Prime Number Theorem and some trivial inequalities,with a bit of calculation that $(r+1)p_1p_2..p_r<n$ for sufficiently large $n$.\hspace{13.5cm}$\Box$   

\vspace{3mm}We don't prove the upper bound here,but we provide a brief sketch.We state a very trivial lemma here.

\vspace{3mm}

\textbf{Lemma 1 :} Given any constant $0<c<\dfrac{6}{\pi^2}$,\hspace{1mm}$\exists$\hspace{1mm} a lattice point $P$ in $A_n$ such that the number 

\hspace{2.25cm}of lattice points from $A_n$ visible from $P$ s at least $cn^2$ if $n$ is sufficiently large.

\vspace{3mm}The proof comes directly from the problem of visibility of lattice points from the origin.We won't dwell with it here.Suppose for $P$ in $A_n$ let $B(P)$ denote the number of points in $A_n$ visible from $P$.So,according to Lemma 5,$|B(P)|>cn^2$,for sufficienty large $n$.Now,we define the points $P_1,P_2...P_k$ recursively as follows : $P_1=(0,0)$ and after the points $P_2,P_3...P_{j-1}$ are chosen,we choose$P_j$ such that $|B(P_j)/\cup_{i=1}^{j-1}B(P_i)|$ is maximum\footnote{For two sets A and B,A/B denotes those elements in A which are not in B}.There may be several choices for $P_j$ but it doesn't matter in this context.The recurrence terminates at the $k$th step where $k$ is the least integer such that $\cup_{i=}^kB(P_i)=A_n$.It is trivial that $A_n$ is visible from the set of points $P_1,P_2...P_k$.We then prove the result by showing that $k\log n$.We will skip this part here.

We use this result of Abbott to prove another beautiful result regarding this topic.

\subsection{On finding an explicit set $B_n$}
\rhead{\textit{Debmalya Basak}}
\lhead{\textit{ON FINDING AN EXPLICIT SET $B_n$}}
\hspace{0.6cm}The next question that obviously came up was finding an explicit set $B_n$.So,let's look at that problem formally.

Let $\omega{(n)}$ denote the number of distinct prime factors of $n,n \in\mathbb{N}$.Suppose $g$($n$) be a function defined on the set of natural numbers such that $g$($n$)$\to \infty$ as $n \to \infty$.We define a set $E_n$($g$)$=$\{$m$ : $m\in\mathbb{N}$,$m \leq n$,$\omega{(m)}\geq g$($n$)\}Now.let's look at the theorem.

\vspace{3mm}
\hspace{1.2cm}\textbf{Theorem 5 :} We can provide an explicit description of a subset $B_n$ of $A_n$ such that $A_n$ 

\hspace{3.75cm}is completely visible from $B_n$ except at most $100|E_n$($g$)$|$ points and for 

\hspace{3.75cm}all sufficiently large $n$,we have $|B_n| \leq 800g$($n$)$\log \log g$($n$).

\vspace{3mm}

The proof is a little intricate.We try to do it step by step.

\vspace{3mm}\hspace{2.35cm}\textit{Proof }: Let $n$ be sufficiently large and let $s=$[10$g$($n$)],$t=$[10$\log \log g$($n$)] and $t_0=$[$\dfrac{t}{10}$]$+1$,where [$x$] denotes the Box function.

Let $X_i=$\{$m$ : $m$ is an integer and $(i-1)t+1<m\leq it$,$i=1,2,3....$\}

\hspace{0.1cm}$I=$\{$i$ : $|X_i\cap E_n$($g$)$| \geq t_0$\} and $Y_i=$\{$it-a$ : $a\in X_i\cap E_n$($g$)\}

Suppose that $1 \leq i \leq n/t$,$i \notin I$ and $1 \leq j \leq n/s$.Then clearly,cardinality of $X_i\cap E_n$($g$) is less than $t_0$.Hence,$|Y_i|<t_0$.Also,suppose $a\notin Y_i$.This means ($it-a$)$\notin X_i\cap E_n$($g$).Let ($it-a$) $\in X_i$.Then ($it-a$)$\notin E_n$($g$).So,$\omega {(it-a)}<g$($n$) if and only if $1\leq a \leq t$ and $a\notin Y_i$.

We now want to count the number of pairs ($a,b$),$1 \leq a \leq t$,$1\leq b \leq s$,$a\notin Y_i$ for which the gcd (($it-a$),($js-b$))$>1$.So,we get,

\vspace{3mm}\hspace{0.6cm}$\sum_{\substack{a=1 \\ a\notin Y_i}}^t\sum_{b=1}^{s}\sum\limits_{\substack{((it-a),(js-b))>1}}1 \hspace{3mm}  \leq \hspace{3mm}\sum_{\substack{a=1 \\ a\notin Y_i}}^t\sum_{b=1}^{s}\sum_p\sum_{\substack{p|(it-a)\\p|(js-b)}}1$\hspace{3mm}where $p$ is a prime.

\vspace{2mm}\hspace{1cm}$\leq$\hspace{3mm}\Bigg($ \sum_{p \leq s}\sum_{\substack{a=1 \\ a\notin Y_i\\p|(it-a)}}^t\sum_{\substack{b=1 \\ p|(js-b)}}^s1$\Bigg)$+$\Bigg($ \sum_{p>s}\sum_{\substack{a=1 \\ a\notin Y_i\\p|(it-a)}}^t1$\Bigg)

\vspace{3mm}
The first inequality comes from the fact that since the gcd is greater than 1,\hspace{1mm}$\exists$\hspace{1mm} at least one prime dividing both ($it-a$) and ($js-b$).Hence the count of this pair goes under the summation of that prime.Now,this count may appear multiple times if the gcd has more than one prime factor.So,the inequality holds.In the 2nd step,we divide the summation by taking primes $\leq s$ on one side and primes greater than $s$ on the other.In the 2nd part of the summation,we are looking for primes greater than $s$ which will divide numbers lying between ($js-1$) and ($j-1$)$s$.So,the count of numbers lying between ($js-1$) and ($j-1$)$s$ is $s$.So,any fixed prime greater than $s$,if chosen,will divide at most one of these $s$ numbers.Hence,we can replace the summation by 1.Proceeding,we get,

\vspace{3mm}\hspace{0.6cm}$\sum_{\substack{a=1 \\ a\notin Y_i}}^t\sum_{b=1}^{s}\sum\limits_{\substack{((it-a),(js-b))>1}}1 \hspace{3mm}  \leq \hspace{3mm} \sum_{p \leq s}$\Big($\dfrac{t}{p}+1$\Big)\Big($\dfrac{s}{p}+1$\Big)$+$\Bigg($ \sum_{p>s}\sum_{\substack{a=1 \\ a\notin Y_i\\p|(it-a)}}^t1$\Bigg)

\vspace{2mm}\hspace{1cm}$\leq$\hspace{3mm}$\sum_{p \leq s}$\Big($\dfrac{t}{p}+1$\Big)\Big($\dfrac{s}{p}+1$\Big)$+$\Bigg($ \sum_{\substack{a=1 \\ a\notin Y_i}}^t1\sum_{p|it-a}$\Bigg)\hspace{3mm} (Removing the restiction on $p$)

\vspace{2mm}\hspace{1cm}$\leq ts\sum_{p}\dfrac{1}{p^2}+(t+s)\sum_{p \leq s}\dfrac{1}{p}+s+(t-|Y_i|)g$($n$)

\vspace{3mm}The foremost inequality is pretty trivial.We are looking at how many numbers lying between ($i-1$)$t$ and ($it-1$) is divisible by $p$.Again,by the same process,there are $t$ such numbers in the interval and $p$ can divide at most [$\dfrac{t}{p}+1$] numbers.We follow the same path for $b$ as well and get that the maximum count is [$\dfrac{s}{p}+1$].So,we replace them in the inequality by \Big($\dfrac{t}{p}+1$\Big)\Big($\dfrac{s}{p}+1$\Big).We just remove the restriction on $p$ in the next step.In the final step,we first take the summation of $\dfrac{1}{p^2}$ for all $p$ and then replace the count of primes $\leq s$ by $s$.On the right most part,we know that ($it-a$)$\notin E_n$($g$),so the  number of distinct primes dividing ($it-a$) is less than $g$($n$).Hence,we get,

\vspace{3mm}\hspace{0.6cm}$\sum_{\substack{a=1 \\ a\notin Y_i}}^t\sum_{b=1}^{s}\sum\limits_{\substack{((it-a),(js-b))>1}}1 \hspace{3mm} < \hspace{3mm} \dfrac{7}{10}ts+(t+s)(\log \log s+O(1))+(t-|Y_i|)g(n)$

\vspace{2mm}\hspace{1cm}$<(t-|Y_i|)s$,for sufficiently large $n$

.\vspace{3mm}The first bound is trivial,taking the sum to be les than $\dfrac{7}{10}$ and we replace the sum of the reciprocals of the primes $\leq s$ by $\log \log s+O(1)$,using a result from Euler's paper on the partial sums of the prime harmonic series,published in 1737.The last step is highly non-trivial.We observe a few properties.
\vspace{3mm}
Given any $m \in \mathbb{N}$,for sufficiently large $n$,we have $t<\dfrac{s}{10^m}$.Also,we have,given any $k \in \mathbb{N}$,for sufficiently large $n$,we have $\log \log s+O(1)< \dfrac{10^kt}{10^{k+1}-1}$.Now,

\vspace{3mm}$\dfrac{7}{10}ts+(t+s)(\log \log s+O(1))+(t-|Y_i|)g(n) \leq \hspace{3mm} \dfrac{7}{10}ts+(t+s)(\log \log s+O(1))+(t-|Y_i|)\dfrac{(s+1)}{10}$

\vspace{2mm}So,we need to show that,

\vspace{3mm}$\dfrac{7}{10}ts+(t+s)(\log \log s+O(1)) < (t-|Y_i|)s-(t-|Y_i|)\dfrac{(s+1)}{10}=(t-|Y_i|)\dfrac{(9s-1)}{10}$

\vspace{3mm}
A few calcuations bring the inequaton down to,

\vspace{3mm}$(t+s)(\log \log s+O(1))+\dfrac{(t-|Y_i|)}{10}<\dfrac{11ts}{100}$

\vspace{3mm}So,it is sufficient to show that,

\vspace{3mm}($s+\dfrac{s}{10^m}$)($\dfrac{10^kt}{10^{k+1}-1}$)$+\dfrac{t}{10}<\dfrac{11ts}{100}$ \hspace{3mm}(As,$|Y_i|\geq 0$)

\vspace{3mm}Simplifying,we have,

\vspace{3mm}($1+\dfrac{1}{10^m}$)($\dfrac{1}{10-\dfrac{1}{10^k}}$)$<\dfrac{11}{100}-\dfrac{1}{10s}$

\vspace{3mm}We can see that as $m,n \to \infty$,the limiting value of the LHS is $\dfrac{1}{10}$.And,we know,for sufficiently large values of $n$,$s>10$.So,the RHS is greater than $\dfrac{1}{10}$ for sufficiently large values of $n$.Hence we can choose $m,n$ large enough so that the ineuality holds for sufficiently large values of $n$.

So,we see that the number of pairs ($a,b$),$1 \leq a \leq t,1 \leq b \leq s$,$a\notin Y_i$ such that (($it-a$),($js-b$))$>1$ is strictly less than $(t-|Y_i|)s$.So,\hspace{1mm}$\exists$\hspace{1mm} integers $a_{ij},b_{ij}$ such that $1 \leq a_{ij} \leq t,1 \leq b_{ij} \leq s$,$a_{ij}\notin Y_i$ and (($it-a_{ij}$),($js-b_{ij}$))$=1$

Now,let's take a lattice point ($x,y$) from $A_n$.We first divide both the X and the Y axis into intervals of length $t$.So,there are non-negative integers $u$ and $v$ such that $ut\leq x <(u+1)t$,$vt\leq y<(v+1)t$.\\

We now consider the 2 following cases.

Case 1 : $u \notin I$.
Let $j$ be the integer such that $js<y\leq (j+1)s$.Now,a lattice point is always visible from any lattice point on it's adjacent stright line.So,f,$u=0$ or $j=0$,then $0 \leq x<t$ or $0 \leq y<s$.So,($x,y$) is visible from the set \{($a,b$) : $1 \leq a \leq 2t,1 \leq b \leq 2s$\}.If $u \geq 1$ and $j \geq 1$,then we can use the result we proved before and say that \hspace{1mm}$\exists$\hspace{1mm} a lattice point ($a_{uj}$,$b_{uj}$) such that $1 \leq a_{uj} \leq t,1 \leq b_{uj} \leq s$,$a_{uj}\notin Y_u$ and (($ut-a_{uj}$),($js-b_{uj}$))$=1$.Hence ($x,y$) is visible from the point ($x-ut+a_{uj}$),($y-js+b_{uj}$).But,$1 \leq x-ut+a_{uj} \leq t+a_{uj} \leq 2t$ and $1 \leq y-js+b_{uj} \leq s+b_{uj} \leq 2s$.So,($x,y$) is visible from the set \{($a,b$) : $1 \leq a \leq 2t,1 \leq b \leq 2s$\}.

Case 2 : $v\notin I$
By symmetry,we follow the same way and say that ($x,y$) is visible from the set \{($a,b$) : $1 \leq a \leq 2s,1 \leq b \leq 2t$\}

So,those points which cannot be seen from \{($a,b$) : $1 \leq a \leq 2t,1 \leq b \leq 2s$\} $\cup$ \{($a,b$) : $1 \leq a \leq 2s,1 \leq b \leq 2t$\}  are among those points for which both $u$ and $v$ belong to $I$.For each fixed ($u,v$),$x$ and $y$ can take t values each independently.So,maximum number of points having same ($u,v$) is $t^2$.So,the number of exceptional points is $t^2|I|^2$ because $u$ and $v$ can each take at most $|I|$ values independently.

Now,another thing to observe is that all $X_i$ are distinct.So,for each entry in $I$,$E_n(g)$ has at least $t_0$ distinct points in it.So,$t_0|I| \leq E_n(g)$.Hence,the number of exceptional points is at most ${\dfrac{t}{t_0}}^2|E_n(g)|^2<100|E_n(g)|^2$.

Now,we need to find the upper bound for the cardinality of $B_n$=\{($a,b$) : $1 \leq a \leq 2t,1 \leq b \leq 2s$\} $\cup$ \{($a,b$) : $1 \leq a \leq 2s,1 \leq b \leq 2t$\}.

Clearly,

\vspace{3mm}\hspace{0.6cm}
$|B_n|<4t^2+2t(2s-2t)+2t(2s-2t)<8st\leq 800g$($n$)$\log \log g$($n$)

\vspace{3mm}Hence,we gave an explicit description of the set $B_n$ with an upper bound of the cardinality of $B_n$.\hspace{15.6cm}$\Box$

\subsection{Corollary of Theorem 5}
\rhead{\textit{Debmalya Basak}}
\lhead{\textit{ON FINDING AN EXPLICIT SET $B_n$}}

Corollary : We can give explicit description of a subset $B_n$ of $A_n$ such that for large $n$,$A_n$ is visible from $B_n$ for at most $100n^2/(\log \log n)^2$ exceptional points,where $|B_n|=O((\log \log n)(\log \log \log \log n))$

\vspace{3mm}\textit{Proof }: We take $g$($n$)$=2\log \log n$.Hardy and Ramanujan proved that $|E_n(g)| \leq \dfrac{n}{\log \log n}$.The rest follows from Theorem 5.\hspace{11.5cm}$\Box$

\vspace{2mm}We state another interesting result here.If we take $g$($n$)$=\dfrac{2(\log n)}{\log \log n}$,$E_n(g)$ becomes empty for sufficiently large $n$.Hence,we can explicitly find a subset $B_n$ such that $A_n$ is completely visible from $B_n$ and $|B_n| \leq 1600\dfrac{(\log n)(\log \log \log n)}{\log \log n}$.

This result is just the application of Theorem 5 except the part that $|E_n(g)|$ becomes empty for large $n$.We try to prove this part.

There is a nice proof of the following inequality that $\omega{(n)}<\dfrac{2(\log n)}{\log \log n}$ for all $n>m$,$n,m \in \mathbb{N}$.We use this inequality here.We look at the finite number of $n$'s for which this inequality doesn't hold and find the maximum value of $\omega{(n)}$ among these.Suppose the maximum value is $\omega{(k)}$ for some $k \in \mathbb{N}$.Since,the function $\dfrac{2(\log n)}{\log \log n}$ is increasing as $n \to \infty$,\hspace{1mm}$\exists$ $l \in \mathbb{N}$ such that $\dfrac{2(\log l)}{\log \log l}>\omega{(k)}$.We want to prove that $\omega{(d)}<g(n)$ for all $d \leq n$,for sufficiently large $n$.We take $n>$max($m,l$).Then for all $d\leq n$,then $\omega{(d)}<g(n)$.This is because $d$ must be either less than $m$ or $l$.If $d<m$,then $\omega{(d)}\leq \omega{(k)}<g(l)<g(n)$.Again,if $d<l$,then $\omega{(d)}<\omega{(l)}<\omega{(n)}$.So,for sufficientl large $n$,$E_n(g)$ becomes empty.\hspace{13.5cm}$\Box$

We end our discussions on this topic here.It is a very interesting topic,with lots of beautiful results and great application of the properties of the Euler Totient function.A lot of credit goes to Dr. Sukumar Das Adhikari for his immense contribution in this field.
\section{Questions that we can look into}
\lhead{\textit{QUESTIONS THAT WE CAN LOOK INTO}}
\rhead{\textit{Debmalya Basak}}
\hspace{0.6cm}Very interesting questions can be asked from all the topics that we have looked into in this article.The concepts can be extended in various ways.Suppose,we proved at the beginning that if $n$ is abundant,then $s$($n$) is abundant for almost all $n$.Then the value of the function $h$($n$) is greater than 2
for both $n$ and $s$($n$).We can study how $h$($n$) and $h$($s$($n$)) varies with respect to $n$.We can look at the difference of these two values and see whether $h$($s$($n$)) is greater than $h$($n$) for almost all $n$ or not.Also,we can work on the range of $\sigma{(n)}$ and $\phi{(n)}$ and see for which values of $n$,$\phi{(n)}$ divides $\sigma{(n)}$ and any special properties of the quotient i.e. can the quotient be a perfect power or not.Also,we can try to study the recurrence reltion for $\sigma{(n)}$ in great details.We can definitely extend the idea of Hidden forests in 3 D.We can also look how the starting point of this Hidden forests vary with $n$ i.e when $n$ is a prime or a highly composite number or a pseudoprime,where does this forest occur at the earliest in the 2 D space.

\section{Bibliography}
\lhead{\textit{ACKNOWLEDGEMENT AND BIBLIOGRAPHY}}
\rhead{\textit{Debmalya Basak}}
1) On a question regarding the visibility of lattice points III,Sukumar Das Adhikari and Yong-

Gao-Chen\\
2) Dynamcs of Sets of zero density,Yuri Gomes Lima\\
3) Asymptotic formula for some functions of Prime numbers,J.Barkley Rosser and Lowell Schoenfield\\
4) An analogue of Grimm's problem of finding distinct prime factors of consecutive integers,Paul 

Erd\H{o}s and Carl Pomerance\\
5) On the product of primes not exceeding $n$,Leo Moser\\
6) On the ratio of the sum of Divisors and Euler Totient function I,Kevin A.Broughan and Daniel 

Delbourgo\\
7) Proof of Euler's Totient Function formula,Shashank Chorge,Juan Vargas\\
8) Summing up the Euler Totient function,Paul Loomis,Michael Plytage and John Polhill\\
9) The Euler Totient,the M\H{o}bius and the divisor functions,Rosica Dineva\\
10) Some results in combinatorical geometry,H.L.Abbott\\

\end{document}